\documentclass[11pt,twoside,english,a4paper]{article}

\usepackage{babel,amssymb,amsthm}
\usepackage{amsmath}
\usepackage{bm}
\usepackage{multirow}
\newtheorem{proposition}{Proposition}[section]

\newtheorem{theorem}[proposition]{Theorem}
\newtheorem{remark}[proposition]{Remark}

\newcommand\smallfrac[2]{\textstyle{\frac{#1}{#2}}}

\begin{document}
\title{A fully discrete Calder\'{o}n Calculus\\
for two dimensional time harmonic waves}

\author{V\'\i ctor Dom\'\i nguez\footnote{Departamento de Ingenier\'\i a
Matem\'atica e Inform\'atica, Universidad P\'ublica de Navarra, 31500 Tudela,
Spain. {\tt victor.dominguez@unavarra.es}. Partially supported by MICINN Project
MTM2010-21037}, Sijiang L. Lu\footnote{Department of Mathematical Sciences,
University of Delaware, USA. {\tt sjly@math.udel.edu}} \& Francisco--Javier
Sayas\footnote{Department of Mathematica Sciences, University of Delaware,
Newark, DE 19716, USA. {\tt fjsayas@math.udel.edu}. Partially supported by NSF
grant DMS 1216356.}}

\maketitle
\centerline{\em Dedicated to Francisco `Paco' Lisbona  on the occasion of his
65th birthday}

\begin{abstract}
In this paper, we present a  fully discretized
Calder\'{o}n Calculus for the two dimensional Helmholtz equation. This full
discretization can be understood as highly non-conforming Petrov-Galerkin
methods, based on two staggered grids of mesh
size $h$, Dirac delta distributions  substituting acoustic charge densities and
piecewise constant functions for approximating acoustic dipole densities. The
resulting numerical schemes from  this calculus are all of order $h^2$ provided
that the continuous equations are well posed. We finish by presenting some
numerical experiments illustrating the performance of this discrete calculus.
\end{abstract}


\section{Introduction}

In this paper we present a very simple and compatible Nystr\"om discretization
of all boundary integral operators for the Helmholtz equation in a smooth
parametrizable curve in the plane. The discretization uses a naif quadrature
method for logarithmic integral equations, based on two staggered grids, and due
to Jukka Saranen and Liisa Schroderus \cite{SaSo:1994} (see also
\cite{SlBu:1992} and \cite{CeDoSa:2002}). This is combined with an equally
simple staggered grid discretization of the hypersingular operator, recently
discovered in \cite{DoLuSa:2012}. If the displaced grids used for the
discretization of these two operators are mutually reversed, then it is possible
to combine these two discretizations with a simple minded Nystr\"om method for
the double layer operator and its adjoint. The complete set of operators is
complemented with a fully discrete version of the single and double layer
potentials. We will explain the construction of the discrete set and reinterpret
it as  a non-conforming Petrov-
Galerkin discretization of the operators (using Dirac deltas and piecewise
constant functions) to which we apply midpoint integration in every element
integral. 

Once the semivariational form has been reached we will show inf-sup conditions
for all discrete operators involved and consistency error estimates based on
asymptotic expansions of the error in the style of \cite{CeDoSa:2002,
DoSa:2001b, DoSa:2001}. We will finally state and sketch the proof of some
convergence error estimates. While some of the results, for individual equations
(mainly based on indirect boundary integral formulations) had already appeared
in previous papers, this is the first time that the entire Calder\'on Calculus
is presented in its entirety. Let it be emphasized, that this is probably the
simplest form of {\em discretizing simultaneously all the potentials and
integral operators} for the Helmholtz equation in the plane and that the methods
we obtain are of order two. Barring the conceptual difficulty of understanding
the boundary integral operators, the methods have the simplicity of basic Finite
Difference Methods and require no effort in their implementation: all discrete
elements are described in full, natural data structures can be easily figured
out from the way the geometry is sampled, and no additional discretization step
(quadrature, assembly by element, mapping to a reference element) is required.
The methods will be presented for the case of a single curve, but we will hint
at its immediate extension to the case of multiple scatterers.

In a final section devoted to numerical experiments, we will show how to use the
methods for transmission problems and how to construct combined field integral
representations.

\section{Calder\'{o}n calculus for exterior Helmholtz boundary problems}

\subsection{Potentials and operators}
Let $\Gamma$ be a smooth simple closed curve given by a regular $1$-periodic
positively oriented parametrization $\mathbf x=(x_1,x_2):\mathbb R \to \Gamma
\subset \mathbb R^2$. Let $\mathbf n(t):=(x_2'(t),-x_1'(t))$ be a non-normalized
outward pointing normal vector at $\mathbf x(t)\in \Gamma$. The domain exterior
to $\Gamma$ will be denoted $\Omega^+$. As a reminder of the fact that we are
taking limits from this exterior domain, the superscript $+$ will be used in
trace and normal derivative operators.

Given $1$-periodic complex-valued functions $\eta$ and $\psi$, the
(parametrized) single and double layer potentials are defined with the formulas
\begin{eqnarray}
\big(\mathrm S\, \eta\big)(\mathbf z)&:=& \frac{\imath}{4}\int_{0}^1
H_0^{(1)}(k|\mathbf z-\mathbf x(t)|)\eta(t)\,\mathrm dt, \label{eq:SL}\\
\big(\mathrm D\, \psi\big)(\mathbf z)&:=& \frac{\imath k}{4}\int_{0}^1
H_1^{(1)}(k|\mathbf z-\mathbf x(t)|)\frac{(\mathbf z-{\bf
x}(t))\cdot\mathbf n(t)}{|\mathbf z-\mathbf n(t)|}\psi(t)\,\mathrm dt \label{eq:DL}
\end{eqnarray}
for arbitrary $\mathbf z\in \mathbb R^2\setminus\Gamma$. (Here $H^{(1)}_n$ is
the Hankel function of the first kind and order $n$.) The single and double
layer potentials define radiating solutions of the Helmholtz equation, namely,
if $U=\mathrm S \eta+\mathrm D \psi$, then
\begin{equation}\label{eq:Helmholtz}
 \Delta U+k^2 U=0  \mbox{ in $\Omega_+$}, \quad \nabla U(\mathbf z)\cdot (\smallfrac1{|\mathbf z|}\mathbf z)
-\imath k U(\mathbf z)=o(\smallfrac1{\sqrt{| \mathbf z |}}),\quad\mbox{as
$|\mathbf z|\to \infty$}.
\end{equation}
Moreover, if $U$ is a $\mathcal C^1(\overline{\Omega^+})$ solution of
\eqref{eq:Helmholtz} and we define 
\begin{equation}\label{eq:cauchyData}
  \varphi= \gamma^+  U:=U|_\Gamma \circ\mathbf x,\qquad 
  \lambda=\partial_\mathbf n^+ U:= ( ( \nabla U)|_\Gamma \circ\mathbf
x)\cdot\mathbf n ,
\end{equation}
then
\cite{HsWe:2010,SaVa:2002}
\begin{equation}
\label{eq:Representation}
 {U}(\mathbf z)= (\mathrm D\,\varphi)(\mathbf z)- 
 (\mathrm S\,\lambda)(\mathbf z), \quad \mathbf z\in\Omega^+. 
\end{equation}
We note that the representation formula \eqref{eq:Representation}, depending on
parametrized Cauchy data \eqref{eq:cauchyData}, can be extended to any locally
$H^1$ solution of \eqref{eq:Helmholtz}. In this work we will restrict our
attention to smooth solutions though.

Associated to the layer potentials we have three integral operators. 
\begin{subequations}
\label{eq:Layers}
\begin{eqnarray} 
({\rm
V}\eta)(s)&:=&\frac{\imath}{4}\int_0^1H_0^{(1)}(k|\mathbf{x}(s)-{
\bf 
x}(t)|)\eta(t)\,\mathrm dt,\label{eq:V}\\ 
(\mathrm K\psi)(s)&:=&  
\frac{\imath k}{4}\int_0^1 
 H_1^{(1)}(k|\mathbf x(s)-\mathbf x(t)|) \frac{(\mathbf x(s)-{\bf 
x}(t))\cdot\mathbf n(t)}{|\mathbf x(s)-\mathbf x(t)|} 
\psi(t)\,\mathrm dt,\label{eq:K}\\ 
(\mathrm{J}\eta)(s)&:=&  
\frac{\imath k}{4}\int_0^1 
 H_1^{(1)}(k|\mathbf x(s)-\mathbf x(t)|) \frac{(\mathbf x(t)-{\bf 
x}(s))\cdot\mathbf n(s)}{|\mathbf x(s)-\mathbf x(t)|} 
\eta(t)\,\mathrm dt, \label{eq:J}
\end{eqnarray}
as well as the integrodifferential operator
\begin{equation}
\label{eq:W}
\mathrm W\psi:= -( \mathrm V \psi')' -k^2 \mathrm V_{\mathbf n}\psi,
\end{equation}
where
\[
(\mathrm V_{\mathbf n}\psi)(s):=\frac{\imath}{4}\int_0^1H_0^{(1)}(k|\mathbf{x}(s)-{\bf 
x}(t)|) \big(\mathbf n(t)\cdot\mathbf n(s)\big)\psi(t)\,\mathrm dt.
\]
\end{subequations}
The operators in \eqref{eq:Layers} are respectively called single layer, double 
layer, adjoint double layer, and hypersingular operator. The operator $\mathrm
W$ admits a different expression in terms of finite parts integrals  (see
\cite[Lemma 2.5.6]{SaVa:2002}), which is where its name comes from.

Layer operators and potentials are related via the so-called  jump
relations \cite{HsWe:2010,McLean:2000,SaVa:2002}, namely, the exterior
parametrized boundary values of the layer operators are given by the formulas
\begin{equation}\label{eq:traces} 
\begin{array}{rclrcl}
\gamma^+ \mathrm S\,\eta&=& \mathrm V \eta,\quad&  
 \gamma^+ \mathrm D\,\psi &=& \smallfrac12\psi+\mathrm K
\psi,\\
 \partial^+_\mathbf n\mathrm S\,\eta&=& -
\smallfrac12 \eta+\mathrm J  \eta,\quad&
 \partial^+_\mathbf n \mathrm D\,\psi&=&-\mathrm W\psi.  
\end{array}
\end{equation}
The matrix of operators
\begin{equation}\label{eq:def:C}
\mathcal{C}^+:=\begin{bmatrix}
           \frac{1}2{\rm I}+\mathrm K& -\mathrm V\\
           -\mathrm W& \frac12{\rm I}-\mathrm J
          \end{bmatrix}
\end{equation}
is the exterior Calder\'on projector. It follows from \eqref{eq:Representation}
and \eqref{eq:traces}, that if  $(\varphi,\lambda)$ are the parametrized Cauchy
data \eqref{eq:cauchyData} for a solution of  \eqref{eq:Helmholtz}, then
$\mathcal C^+(\varphi,\lambda)^\top=(\varphi,\lambda)^\top$ or, equivalently
\begin{equation}
\label{eq:C:calculus}
\mathcal D^+\begin{bmatrix}\varphi \\ \lambda\end{bmatrix}:=
\begin{bmatrix}
           \frac12\mathrm I-\mathrm K& \mathrm V\\
           \mathrm W& \frac12\mathrm I+\mathrm J
 \end{bmatrix}
\begin{bmatrix}\varphi \\ \lambda\end{bmatrix}=\begin{bmatrix} 0 \\ 0 \end{bmatrix}.
\end{equation}
 Note that $\mathrm K$ and $\mathrm J$ are transposed of each other, while $\mathrm V$ and $\mathrm W$ are symmetric.

\subsection{Boundary integral equations for exterior problems}
We next summarize a collection of boundary integral equations leading to the solution of  \eqref{eq:Helmholtz}
with a given boundary condition:
\begin{equation}\label{eq:BC}
\gamma^+ U=\beta_0 \qquad \mbox{or}\qquad \partial_{\mathbf n}^+ U=\beta_1.
\end{equation}
The data functions in the right-hand side of \eqref{eq:BC} are 
$1$-periodic functions and the boundary operators are those of
\eqref{eq:cauchyData}. Recall that the Dirichet or Neumann exterior
problem  for the Helmholtz equation with Sommerfeld radiation condition
at infinity are uniquely solvable. 

A {\em direct method} for solving the exterior Dirichlet problem starts
in the representation formula \eqref{eq:Representation}, equates
$\varphi=\beta_0$, and then uses one of the two identities in 
\eqref{eq:C:calculus} to set up an integral equation in order to find $\lambda$.
Similarly, for the Neumann problem, we impose  $\lambda=\beta_1$, and then use
one of the equations in \eqref{eq:C:calculus} in search of $\varphi$. The
resulting integral equations are collected in Table~\ref{tab:d}.

\begin{table}[ht]
\[
 \begin{array}{|l|lll|}
  \hline 
  \multirow{5}{*}{\text{Dirichlet}} &  & & \ \\
  &\mathrm V \lambda=  -\smallfrac12\varphi+\mathrm K \varphi& \varphi=\beta_0 &
\text{(dD01)}
 \\[1.25ex]
 \cline{2-4} & & &\\
& \smallfrac12\lambda+\mathrm J\lambda=-\mathrm W\varphi, &\varphi=\beta_0 &\text{(dD02)}
\\[1.25ex]
\hline
  \multirow{5}{*}{\mbox{Neumann}} &   & & \ \\
  &-\smallfrac12\varphi+\mathrm K \varphi=\mathrm V\lambda, &  \lambda=\beta_1 & \text{(dN01)} \\[1.25ex]
 \cline{2-4} && & \\
&-\mathrm W\varphi= \smallfrac12\lambda+\mathrm J \lambda,& \lambda=\beta_1 & \text{(dN02)}\\[1.25ex]
\hline 
 \end{array}
\]
\caption{BIEs for direct formulations. The representation formula
is \eqref{eq:Representation}. All these equations are solvable. Uniqueness is
discussed in Proposition \ref{prop:uniq}.}\label{tab:d}
\end{table}

An {\em indirect method} based on the single layer potential representation
looks for $U=\mathrm S \eta$ and then uses the expressions in the first column
of \eqref{eq:traces} to set up an integral equation depending on which boundary
data is known. Similarly, we can look for $U=\mathrm D\psi$ and use the boundary
integral operators that appear in the right column of \eqref{eq:traces} to build
an equation. These equations are gathered in Table~\ref{tab:i}.

 \begin{table}[ht]
\[
 \begin{array}{|l|lll|}
  \hline 
\multirow{5}{*}{\text{Dirichlet}} &  & & \ \\
  &\mathrm V \eta=  \beta_0, &  U=\mathrm S\,\eta&
\text{(iD01)}  \\[1.25ex]
 \cline{2-4} & & &\\
& \smallfrac12\psi+\mathrm K\psi=\beta_0,& U=\mathrm D\,\psi 
&\text{(iD02)} \\[1.25ex]
  \hline
\multirow{5}{*}{\text{Neumann}} &  & & \ \\
  &-\smallfrac12\eta+\mathrm J \eta =  \beta_1,& U=\mathrm S\,\eta&
\text{(iN01)}
 \\[1.25ex]
 \cline{2-4} && &\\
&\mathrm W\psi =  -\beta_1,& U=\mathrm D\,\psi &\text{(iN02)}\\[1.25ex]  
\hline 
 \end{array}
\]
\caption{BIEs for indirect formulations. The potential representation is given
next to the boundary integral equation. Unique solvability of these equations is
discussed in Proposition \ref{prop:uniq}.}\label{tab:i}
\end{table}

\begin{proposition}[See {\cite[Section 3.2]{ne:2001}}]\label{prop:uniq} 
Let $\Omega$ be the domain interior to $\Gamma$. 
\begin{itemize}
\item[{\rm (a)}] Equations  {\rm (dN01)}, {\rm (iN01)},  {\rm (dD01)}, and {\rm
(iD01)} are uniquely solvable if and only if $-k^2$ is not a Dirichlet
eigenvalue of the Laplace operator in $\Omega$.
\item[{\rm (b)}] Equations  {\rm (dN02)}, {\rm (iN02)},  {\rm (dD02)}, and {\rm
(iD02)} are uniquely solvable if and only if $-k^2$ is not a Neumann eigenvalue
of the Laplace operator in $\Omega$.
\end{itemize}
\end{proposition}

The equations of Tables \ref{tab:d} and \ref{tab:i} involve the four operators
of the matrix $\mathcal D^+$ in \eqref{eq:C:calculus} and their transposes. The
operators in the first row of $\mathcal D^+$ are invertible when $-k^2$ is not
an interior Dirichlet eigenvalue. The operators in the second row of $\mathcal
D^+$ are invertible when $-k^2$ is not an interior Neumann eigenvalue. The
precise Sobolev space setting where these equations are well posed will be
explained in Section \ref{sec:4.1}. In addition to these equations, the
Calder\'on Calculus, given by the jump relations \eqref{eq:traces} and the
identities \eqref{eq:C:calculus} can be used to construct combined integral
equations and several other associated boundary integral equations, some of
which are invertible for all values of $k$. 

\section{The fully discrete calculus}

\subsection{Matrix representation}\label{sec:3}
Let $N$ be a positive integer, $h:=1/N$, and let us consider the uniform grid in
parametric space
\[ 
s_i:=(i-{\textstyle\frac12})h, \quad
 t_i:=i h,\qquad i\in\mathbb{Z},
\]
thus defined so that $t_i$ is the midpoint of the interval $(s_i,s_{i+1})$.  The
following quantities will be {\em all the geometric elements of $\Gamma$} that
will be used in the discrete Calculus: 
\begin{equation}\label{eq:geometry}
\begin{array}{rclrcl}
 {\bf m}_i&:=&\mathbf x(t_i),\quad &{\bf b}_i&:=&\mathbf x(s_i),\\
\mathbf n_i&:=&h \mathbf n(t_i), \quad &{\bf \ell}_i&:=&|\mathbf n_i|=h |{\bf
x}'(t_i)|,
\qquad {\bf s}_i=h^2\mathbf x''(t_i).           
\end{array}
\end{equation}
These quantities make up the {\em main discretization grid}. Note that they are
defined for $i\in \mathbb Z$, modulo $N$. For practical reasons, we will need a
discrete function $n(i)$, that gives the next index in a rotating (modulo $N$)
form, so that $n(i)=i+1$ for $i\le N-1$ and $n(N)=1$.
We now take
$\varepsilon\in(-1/2,1/2)\setminus\{0\}$ and repeat the same construction
with the displaced grid in parametric space:
\[
 t_i^\varepsilon=(i+\varepsilon) h,\qquad
s_i^\varepsilon=(i+\varepsilon-{\textstyle\frac12})h.
\]
The quantities
 ${\bf m}_i^\varepsilon$, ${\bf b}_i^\varepsilon$,  ${\bf
n}_i^\varepsilon$, ${\bf \ell}^\varepsilon_i$,  and
${\bf s}^\varepsilon_i$  are defined accordingly. They constitute the {\em companion grid}.

Given column vectors $\boldsymbol\eta=(\eta_1,\ldots,\eta_N)^\top\in \mathbb
C^N$, $\boldsymbol\psi=(\psi_1,\ldots,\psi_N)^\top \in \mathbb C^N$, we consider
the discrete single and double layer potentials:
\begin{subequations}\label{eq:poth}
\begin{eqnarray}
\mathrm S_h (\mathbf z)\, \boldsymbol\eta &:=& \sum_{j=1}^N\frac{\imath}{4} 
H_0^{(1)}(k|\mathbf z-{\bf m}_j^{\varepsilon}|)\eta_j, \\
\mathrm D_h (\mathbf z) \,\boldsymbol\psi&:=&
\sum_{j=1}^N \frac{\imath k}{4} 
H_1^{(1)}(k|\mathbf z-{\bf m}_j|)\frac{(\mathbf z-{\bf m}_j)\cdot\mathbf n_j}{|{\bf
z}-{\bf m}_j|}\psi_j.
\end{eqnarray}
\end{subequations}
We also consider four $N\times N$  matrices $\mathrm V_h$, $\mathrm K_h$,
$\mathrm J_h$ and $\mathrm W_h$, given by
\begin{subequations} \ \label{eq:discreteOperators}
\begin{eqnarray}
\mathrm V_{ij}&=& \frac{\imath}{4} H_0^{(1)}(k|{\bf m}_i-
{\bf m}_j^{\varepsilon}|),\\ 
 \mathrm{K}_{ij}&:=&\left\{\begin{array}{ll}
            \displaystyle                \frac{{\bf s}_i\cdot{\bf
n}_i}{4\pi\ell_i^2},&i=j,\\
\displaystyle\frac{\imath k}{4} 
 H_1^{(1)}(k|{\bf m}_i-{\bf m}_j|) \frac{({\bf m}_i-
 {\bf m}_j)\cdot \mathbf n_i}{|{\bf m}_i-{\bf m}_j|},\quad&i\neq j,
                           \end{array}\right.\\
 \mathrm{J}_{ij} &:=&  
 \left\{\begin{array}{ll}
 \displaystyle \frac{{\bf s}^\varepsilon_i\cdot{\bf
n}^\varepsilon_i}{4\pi(\ell_i^\varepsilon)^2},&i=j,\\
\displaystyle \frac{\imath k}{4} 
 H_1^{(1)}(k|{\bf m}_i^\varepsilon-{\bf m}_j^\varepsilon|)
\frac{({\bf m}_j^\varepsilon-
{\bf m}_i^\varepsilon)\cdot \mathbf n_j^\varepsilon}
{|{\bf m}_i^\varepsilon-{\bf m}_j^\varepsilon|},\quad&i\neq j,
\end{array}\right.
\\ 
 \mathrm{W}_{ij}&:=&
\widetilde{\mathrm V}_{n(i),n(j)}+\widetilde{\mathrm V}_{ij}-\widetilde{\mathrm V}_{n(i),j}
-\widetilde{\mathrm V}_{i,n(j)}  -k^2 (\mathbf n_{i}^\varepsilon\cdot {\bf
n}_{j}) \mathrm V_{ji},\label{eq:Wh}
\end{eqnarray}
\end{subequations}
where
\[
 \widetilde{\mathrm V}_{ij}=
\frac{\imath}{4} H_0^{(1)}(k|\mathbf{b}_i^\varepsilon-{\bf 
b}_j|). 
\]
Note that the diagonal values in $\mathrm K_h$ and $\mathrm J_h$ are defined
using the limit values in the kernels of the integral operators $\mathrm K$ and
$\mathrm J$ as $|s-t|\to 0$.

\begin{remark}
As can be seen from \eqref{eq:poth} and \eqref{eq:discreteOperators}, the
structure of the matrices and operators does not remember where the discrete
geometric data come from. The formulas \eqref{eq:poth} and
\eqref{eq:discreteOperators} use discrete data $\{\mathbf m_i,\mathbf
n_i,\mathbf b_i, \ell_i, \mathbf s_i\}$ and $\{\mathbf m_i^\varepsilon,\mathbf
n_i^\varepsilon,\mathbf b_i^\varepsilon, \ell_i^\varepsilon, \mathbf
s_i^\varepsilon\}$, sampled from the curve. It is immaterial whether these data
have been sampled from a simple curve or several simple non-intersecting curves.
The next-index function $n(i)$ used in $\mathrm W_h$ has to be adapted to
contain cycles of nodes showing the different connected components of the
collection of curves.
\end{remark}

Discretization of the integral equations in Tables \ref{tab:d} and \ref{tab:i}
is almost straightforward based on these matrices and potentials. The Dirichlet
and Neumann data in \eqref{eq:BC} are discretized by vectors of samples:
\begin{equation}\label{eq:beta}
\boldsymbol\beta_0:=(\beta_0(t_1),\ldots,\beta_0(t_N))^\top \qquad
\boldsymbol\beta_1:=h\,(\beta_1(t_1^\varepsilon),\ldots,
\beta_1(t_N^\varepsilon))^\top.
\end{equation}
The different scaling of these vectors will be clear from the interpretation of
these methods that we will give in Section \ref{sec:3.2}. At this stage, it can
be justified with some arguments of dimensional analysis, given the fact that
$\beta_1$ corresponds to data of a derivative of the function. By the definition
of the parametrized boundary operators \eqref{eq:cauchyData}, of the Cauchy data
\eqref{eq:BC} and of the discrete quantities \eqref{eq:geometry}, we can
similarly write
\[
\boldsymbol\beta_0=(U(\mathbf m_1),\ldots,U(\mathbf m_N))^\top \qquad
\boldsymbol\beta_1=(\nabla U(\mathbf m_1^\varepsilon)\cdot\mathbf
n_1^\varepsilon,\ldots,\nabla U(\mathbf m_N^\varepsilon)\cdot\mathbf
n_N^\varepsilon)^\top.
\]
The {\em discrete direct methods} use a representation formula
\begin{equation}\label{eq:reprh}
U_h(\mathbf z)=\mathrm D_h(\mathbf z)\boldsymbol\lambda-\mathrm S_h(\mathbf z)\boldsymbol\varphi
\end{equation}
and one of the linear systems of Table \ref{tab:dd}. The {\em discrete indirect
methods} appear collected in Table \ref{tab:ii}, including the corresponding
potential representation.

\begin{table}[ht]
\[
 \begin{array}{|l|lll|}
  \hline 
  \multirow{5}{*}{\text{Dirichlet}} &  & & \ \\
  &\mathrm V_h \boldsymbol\lambda=  -\smallfrac12\boldsymbol\varphi+\mathrm K_h \boldsymbol\varphi& \boldsymbol\varphi=\boldsymbol\beta_0 &
\text{(dD01h)}
 \\[1.25ex]
 \cline{2-4} & & &\\
& \smallfrac12\boldsymbol\lambda+\mathrm J_h\boldsymbol\lambda=-\mathrm W_h\boldsymbol\varphi, &\boldsymbol\varphi=\boldsymbol\beta_0 &\text{(dD02h)}
\\[1.25ex]
\hline
  \multirow{5}{*}{\mbox{Neumann}} &   & & \ \\
  &-\smallfrac12\boldsymbol\varphi+\mathrm K_h \boldsymbol\varphi=\mathrm
V_h\boldsymbol\lambda, &
  \boldsymbol\lambda=\boldsymbol\beta_1 & \text{(dN01h)} \\[1.25ex]
 \cline{2-4} && & \\
&-\mathrm W_h\boldsymbol\varphi= \smallfrac12\boldsymbol\lambda+\mathrm J_h \boldsymbol\lambda,& \boldsymbol\lambda=\boldsymbol\beta_1 & \text{(dN02h)}\\[1.25ex]
\hline 
 \end{array}
\]
\caption{Discrete direct methods, with  representation formula \eqref{eq:reprh}.}\label{tab:dd}
\end{table}

\begin{table}[ht]
\[
 \begin{array}{|l|lll|}
  \hline 
\multirow{5}{*}{\text{Dirichlet}} &  & & \ \\
  &\mathrm V_h \boldsymbol\eta=  \boldsymbol\beta_0, &  U_h=\mathrm S_h\,\boldsymbol\eta&
\text{(iD01h)}  \\[1.25ex]
 \cline{2-4} & & &\\
& \smallfrac12\boldsymbol\psi+\mathrm K_h\boldsymbol\psi=\boldsymbol\beta_0,& U_h=\mathrm D_h\,\boldsymbol\psi 
&\text{(iD02h)} \\[1.25ex]
  \hline
\multirow{5}{*}{\text{Neumann}} &  & & \ \\
  &-\smallfrac12\boldsymbol\eta+\mathrm J_h \boldsymbol\eta = \boldsymbol \beta_1,& U_h=\mathrm S_h\,\boldsymbol\eta&
\text{(iN01h)}
 \\[1.25ex]
 \cline{2-4} && &\\
&\mathrm W_h\boldsymbol\psi =  -\boldsymbol\beta_1,& U_h=\mathrm D_h\,\boldsymbol\psi &\text{(iN02h)}\\[1.25ex]  
\hline 
 \end{array}
\]
\caption{Discrete indirect methods.}\label{tab:ii}
\end{table}

\subsection{Reinterpretation as non-conforming Petrov-Galerkin methods}\label{sec:3.2}

Our method can be understood as a collection  of  non-conforming
Petrov-Galerkin methods with a very simple quadrature rule for approximating any
integral appearing in the scheme. The basic idea is the following: the input of
$\mathrm D$ (and therefore $\mathrm W$ and $\mathrm K$) will be approximated
with a piecewise constant function on the main grid; the input of $\mathrm S$
(and therefore $\mathrm V$ and $\mathrm J$) will be approximated with a linear
combination of Dirac deltas on the companion grid; tests related to Dirichlet
problems will be carried out by Dirac deltas on the main grid; test related to
Neumann problems will be done with piecewise constants on the companion grid;
finally, all integrals will be broken into subintervals of the grid and
approximated with a midpoint rule.

In order to write the methods of Section \ref{sec:3} in the form where we will
develop their convergence analysis, we need to define some new discrete
elements. First of all, we consider the (periodic) Dirac delta distribution
$\delta_z$ at a point $z$. Its action on any periodic function that is
continuous around $z$ will be denoted
$\{\delta_z,\rho\}=\{\rho,\delta_z\}:=\rho(z)$. Given an interval $I\subset
\mathbb R$, we will denote by $\chi_I$ the periodized characteristic function of
$I$, i.e., the characteristic function of the set $I+\mathbb Z$. We then
consider four discrete spaces
\[
\begin{array}{rclrcl}
\!\!S_h &\!\!\!:=\!\!\!& \mathrm{span}\{ \chi_{(s_{i-1},s_i)}\,:\,
i=1,\ldots,N\},  &
S_{h,\varepsilon} &\!\!\!:=\!\!\!& \mathrm{span}\{
\chi_{(s_{i-1}^\varepsilon,s_i^\varepsilon)}\,:\, i=1,\ldots,N\},\\
\!\!S_h^{-1} &\!\!\!:=\!\!\!& \mathrm{span}\{ \delta_{t_i}\,:\,
i=1,\ldots,N\},  &S_{h,\varepsilon}^{-1} &\!\!\!:=\!\!\!& \mathrm{span}\{
\delta_{t_i^\varepsilon}\,:\, i=1,\ldots,N\}.\\
\end{array}
\]
For elements of these spaces we will identify the vector of their coefficients
--with respect to the basis that has been used to define the space--, using the
same letter in boldface font. For example,
\[
S_h^{-1} \ni \mu_h=\sum_{j=1}^N \mu_j \delta_{t_j} \longleftrightarrow \boldsymbol\mu=(\mu_1,\ldots,\mu_N)^\top \in \mathbb C^N. 
\]
The two discrete operators 
\[
Q_h^{-1} \rho:=
h \sum_{j=1}^N \rho(t_j)\delta_{t_j} \qquad
Q_{h,\varepsilon}^{-1}\rho:=h\sum_{j=1}^N
\rho(t_j^\varepsilon)\delta_{t_j^\varepsilon} 
\]
complete the collection of elements needed for a more variational description of
the discrete Calder\'on Calculus. They will be used to denote midpoint
quadrature approximations. For example,
\[
\{ Q_h^{-1} \rho,\phi\}=h \sum_{j=1}^N \rho(t_j)\phi(t_j) \approx \int_0^1 \rho(t)\phi(t)\mathrm d t.
\]
The discrete potentials \eqref{eq:poth} can be easily described in this language:
\[
S_{h,\varepsilon}^{-1}\ni \eta_h \mapsto \mathrm S_h(\cdot) \boldsymbol\eta =
\mathrm S \eta_h, \qquad S_h\ni \psi_h \mapsto \mathrm
D_h(\cdot)\boldsymbol\psi=\mathrm D Q_h^{-1}\psi_h.
\]
Observe how in the double layer potential we are just applying the midpoint rule
to approximate $\mathrm D \psi_h$, while no additional integration is needed
in the already fully discrete expression for $\mathrm S \eta_h$.

The matrices \eqref{eq:discreteOperators} have their variational counterparts as
bilinear forms:
\begin{eqnarray*}
 S_h^{-1}\times S_{h,\varepsilon}^{-1} \ni (\mu_h,\eta_h) & \longmapsto & \mathrm v(\mu_h,\eta_h):=\{\mu_h,\mathrm V\eta_h\} =\boldsymbol\mu^\top \mathrm V_h\boldsymbol\eta,\\
 S_h^{-1}\times S_h \ni (\mu_h,\psi_h) & \longmapsto & \mathrm k(\mu_h,\psi_h):=\{\mu_h,\mathrm KQ_h^{-1}\psi_h\} =\boldsymbol\mu^\top \mathrm K_h\boldsymbol\psi,\\
S_{h,\varepsilon}\times S_{h,\varepsilon}^{-1} \ni (\phi_h,\eta_h) &
\longmapsto & \mathrm j(\phi_h,\eta_h):=\{Q_{h,\varepsilon}^{-1}\phi_h,\mathrm
J\eta_h\} =\boldsymbol\phi^\top \mathrm J_h\boldsymbol\eta,\\
S_{h,\varepsilon}\times S_h \ni (\phi_h,\psi_h) &\longrightarrow & \mathrm w(\phi_h,\psi_h):=\boldsymbol\phi^\top \mathrm W_h \boldsymbol\psi.
\end{eqnarray*}
The bilinear form $\mathrm w$ can be understood as follows
\[
\mathrm w(\phi_h,\psi_h)=\{ \phi_h',\mathrm V \psi_h'\} -k^2\{Q_{h,\varepsilon}^{-1} \phi_h,\mathrm V_{\mathbf n} Q_h^{-1}\psi_h\},
\]
just by noticing that $\chi_{(s_{i-1},s_i)}'=\delta_{s_{i-1}}-\delta_{s_i}$ and
that a change of sign has to be applied to the leading integrodifferential part
of $\mathrm W$ (see \eqref{eq:W}) when changing the differentiation to the test
function. The rationale behind this choice of spaces can be observed in the
matrix of operators $\mathcal D^+$ in \eqref{eq:C:calculus}. As trial spaces we
are considering $S_h\times S_{h,\varepsilon}^{-1}$, while the rows of $\mathcal
D^+$ are respectively tested with $S_h^{-1}$ and $S_{h,\varepsilon}$. This means
that the operators of the second kind ($\pm\frac12\mathrm I+\mathrm K$ and
$\pm\frac12\mathrm I+\mathrm J$) are discretized on a single grid (each of them
on a different grid though), while the operators of the first kind ($\mathrm V$
and $\mathrm W$) use two grids. This is actually a requirement due to the fact
that the kernels of $\mathrm V$ and $\mathrm V_{\mathbf n}$ cannot be evaluated
in the diagonal $s=t$, where they have a logarithmic singularity. Once 
this choice of trial and test spaces has been taken as a first step in the
discretization of the four operators in \eqref{eq:C:calculus}, midpoint
integration is applied to all remaining integrals. The operators $Q_h^{-1}$ and
$Q_{h,\varepsilon}^{-1}$ are used as a way of enforcing full discretization of
every operator acting on a piecewise constant function.

To describe variationally the equations in Tables \ref{tab:dd} and \ref{tab:ii}
we first cast the data function ($\beta_0$ for the Dirichlet problem and
$\beta_1$ for the Neumann problem) in the discrete spaces
\[
\beta_0^h:=\sum_{j=1}^N \beta_0(t_j)\chi_{(s_{j-1},s_j)}\in S_h,\qquad
\beta_1^h:=Q_{h,\varepsilon}^{-1}\beta_1=h\sum_{j=1}^N
\beta_1(t_j^\varepsilon)\delta_{t_j^\varepsilon}\in S_{h,\varepsilon}^{-1},
\]
so that their coefficients coincide with the sample vectors \eqref{eq:beta}. The
equations (dN01h) correspond then to writing $\lambda_h=\beta_1^h$, solving
\begin{equation}\label{eq:dN01h}
\varphi_h\in S_h
\quad \mbox{s.t.}\quad -\smallfrac12\{\mu_h,\varphi_h\}+\mathrm
k(\mu_h,\varphi_h)=\mathrm v(\mu_h,\lambda_h)\quad \forall \mu_h \in S_h^{-1},
\end{equation}
and finally using $U_h=\mathrm D Q_h^{-1}\varphi_h-\mathrm S \lambda_h$ as
discrete representation formula.
The indirect method (iN02h) corresponds to solving
\[
\psi_h \in S_h \quad\mbox{s.t}\quad 
\mathrm
w(\phi_h,\psi_h)=-\{\beta_1^h,\phi_h\}=-\{Q_{h,\varepsilon}^{-1}\phi_h,\beta_1\}
\quad \forall \phi_h\in S_{h,\varepsilon},
\]
for a potential representation $U_h=\mathrm D Q_h^{-1} \psi_h$. The indirect method (iD01h) is equivalent to solving
\[
\eta_h \in S_{h,\varepsilon}^{-1} \quad\mbox{s.t.}\quad
\mathrm v(\mu_h,\eta_h)=\{\mu_h,\beta_0^h\} = \{\mu_h,\beta_0\} \quad \forall
\mu_h \in S_h^{-1}.
\]
The remaining five discrete equations in Tables \ref{tab:dd} and \ref{tab:ii}
can be easily rewritten using these same elements.

\section{Numerical analysis}

\subsection{Stability}\label{sec:4.1}

Analysis of the methods in Section \ref{sec:3} is carried out in the form given
in Section \ref{sec:3.2}, in the frame of periodic Sobolev spaces. For $s\in
\mathbb R$ we define the space $H^s$ as the completion of the space of
trigonometric polynomials $\mathrm{span}\,\{ \exp(2\pi\imath
m\,\cdot\,)\,:\,m\in \mathbb Z\}$ with respect to the norm
\[
 \|\rho\|_{s}^2=|\widehat{\rho}(0)|^2+\sum_{m\ne
0}|m|^{2s}|\widehat{\rho}(m)|^2,\quad 
 \widehat{\rho}(m):=\int_0^1 \rho(t)\exp(-2\pi\imath m t)\,\mathrm dt.
\]
An extensive treatment of these spaces can be found in \cite{SaVa:2002}. The
operators \eqref{eq:Layers} can be extended to act on all Sobolev spaces $H^s$.
In particular, the following result holds (see \cite[Table 2.1.1]{HsWe:2010} and
\cite[Section 3.2]{ne:2001}).

\begin{proposition}
The operators
\begin{equation}\label{eq:allOperators}
 \pm {\textstyle\frac 12}+\mathrm K, \ 
 \pm {\textstyle \frac 12}+\mathrm J: H^s\to H^s,\quad  \mathrm V:H^s\to
H^{s+1}, \quad \mathrm W:H^s\to H^{s-1}
\end{equation}
are bounded for all $s$. If, in addition, $-k^2$ is neither a Dirichlet nor a
Neumann eigenvalue of the Laplacian in $\Omega$ (cf. Proposition
\ref{prop:uniq}), then all of them are invertible.
\end{proposition}

\begin{proposition}\label{prop:infsup} Assume that $-k^2$ is neither a Dirichlet
nor a Neumann eigenvalue of the Laplacian in $\Omega$ and let
$\varepsilon\in(-1/2,1/2) \setminus\{0\}$. Then
there exist positive numbers $c_{\mathrm{V}}, c_{\mathrm{K}},c_{\mathrm{J}}, c_{\mathrm{W}}>0$
so that for all $h$ small enough 
\begin{eqnarray}
\label{eq:stabV}
 \inf_ {0\ne \eta_h\in S_{h,\varepsilon}^{-1}
}\sup_{0\ne \mu_h\in S_{h}^{-1}}\frac{|\mathrm v(\mu_h,\eta_h)|}{\|\mu_h\|_{-1}\|\eta_h\|_{-1}}&\ge&
c_{\mathrm{V}},\\
\inf_{0\ne \psi_h\in S_{h}}
\sup_{0\ne \mu_h\in S_{h}^{-1}
}
\frac{| \pm \smallfrac12\{\mu_h,\psi_h\} +\mathrm k(\mu_h,
\psi_h )|}{\|\mu_h\|_{-1}\|\psi_h\|_{0}}\label{eq:stabK}
&\ge&
c_{\mathrm{K}},\\ 
 \inf_{0\ne \eta_h\in S_{h,\varepsilon}^{-1}}
 \sup_{0\ne \phi_h\in
S_{h,\varepsilon}}
\frac{|\pm \smallfrac12\{ \phi_h,\eta_h\} + \mathrm j(
\phi_h,\eta_h)|}{\|\phi_h\|_{0}\|\eta_h\|_{-1}}
&\ge&
c_{\mathrm{J}},\label{eq:stabJ}\\ 
\inf_{0\ne \psi_h\in S_h}
\sup_{0\ne \phi_h\in
S_{h,\varepsilon} }
 \frac{|\mathrm{w}(\phi_h, \psi_h)|}{\|\phi_h\|_{0}\|\psi_h\|_{0}}&\ge&
c_{\mathrm{W}}.\label{eq:stabW}
\end{eqnarray}
The constants can depend on $\varepsilon$.
\end{proposition}
\begin{proof}
Condition \eqref{eq:stabV} was proved in \cite[Proposition 8]{CeDoSa:2002},
although it is based on a stability result (phrased in different terms) given in
\cite{SaSo:1994}.
Condition \eqref{eq:stabW} has been proven in \cite[Theorem 1]{DoLuSa:2012}.
With minor modifications, the proof of \cite[Theorem 2]{DoRaSa:2008} can be used
to prove \eqref{eq:stabJ}. It is then easy to note that this result would also
hold for the spaces $S_h^{-1}$ and $S_h$ (it all amounts to displacing the grid
for both test and trial functions). Then, by an easy transposition argument,
\eqref{eq:stabK} holds.
\end{proof}

The value $\varepsilon=0$ is not a practicable option for the choice of the
grids: in this case both grids coincide and we are obliged to evaluate the
singular kernels in their diagonal. The choices $\varepsilon=\pm1/2$ lead to a
discretization of $\mathrm V$ (they give the same one)  that is not stable,
i.e., the inf-sup condition does not hold. The proof of the inf-sup condition
for the discretization of $\mathrm W$ in \cite{DoLuSa:2012} requires also that
$\varepsilon\neq \pm1/2$, because it is based on the result for $\mathrm V$,
although numerical evidence points to this being just a technical restriction,
which is not in the case of $\mathrm V$. Note finally that dependence of the
methods on $\varepsilon$ is $1$-periodic.

\subsection{Consistency analysis via asymptotic expansions}

We next study the consistency of the approximation of the bilinear forms
associated to the four operators \eqref{eq:Layers} by their discrete
counterparts, as well as the approximation of the identity operators that appear
in the equations of Tables \ref{tab:d} and \ref{tab:i}. The consistency error
analysis is carried out by comparison with a {\em quasioptimal projection} of
the corresponding unknown (the input of the integral operator) in the discrete
space. These projections are defined by matching the central Fourier
coefficients:
\begin{eqnarray*}
 S_{h,\varepsilon}^{-1}&\ni& D_{h,\varepsilon}^{-1}\eta, \quad
\widehat{D_{h,\varepsilon}^{-1}\eta}(m)=\widehat{\eta}(m), \quad
-N/2<m\le N/2, \\
 S_{h}&\ni& D_{h}\psi, \quad
\widehat{D_{h}\psi}(m)=\widehat{\psi}(m), \quad
-N/2<m\le N/2 . 
\end{eqnarray*}
The operator $D_h$ was studied in \cite{arn:1983}, while
$D_{h,\varepsilon}^{-1}$ proceeds from \cite{CeDoSa:2002}. It is proved in those
references that
\begin{subequations}\label{eq:Dh}
\begin{eqnarray}
\label{eq:Dha}
 \|D_h\psi-\psi\|_{s}\le C_{s,r}h^{r-s}\|\psi\|_{r} & & s\le r\le 1,\ s<1/2,\\
\label{eq:Dhb}
\| D_{h,\varepsilon}^{-1}\eta-\eta\|_s \le C h^{r-s}\|\eta\|_r & & s\le r\le 0,\
s<-1/2.
\end{eqnarray}
\end{subequations}

\begin{proposition}\label{prop:exp1}
For all $\eta\in H^3$ and $\psi\in H^4$ it
holds 
\begin{eqnarray*}
|\{\phi_h,D_{h,\varepsilon}^{-1}\eta\} -\{\phi_h,Q_{h,\varepsilon}^{-1}\eta\}
|&\le&Ch^{3}\|\eta\|_3\|\phi_h\|_0,\quad \forall \phi_h\in
S_{h,\varepsilon},\\
|\{\mu_h, D_{h}\psi\}-\{\mu_h,\psi \}+ h^2\smallfrac1{24}\{\mu_h,\psi''\}
|&\le&C h^3\|\psi\|_4\|\mu_h\|_{-1},\quad \forall \mu_h\in
S_{h}^{-1}. 
\end{eqnarray*}
The constants in the bounds are independent of $\varepsilon$.
\end{proposition}
\begin{proof} The second expansion follows from \cite[Theorem 7]{CeDoSa:2002}.
To prove the first one, note that by \cite[Lemma 5]{CeDoSa:2002}
\begin{equation}\label{eq:Eh}
D_{h,\varepsilon}^{-1} \eta-Q_{h,\varepsilon}^{-1}\eta=Q_{h,\varepsilon}^{-1}
E_h\eta, \mbox{ where } E_h\eta:=\sum_{-\frac{N}2<m\le
\frac{N}2}\widehat\eta(m)\exp(2\pi\imath m\cdot)-\eta.
\end{equation}
A direct computation (see also \cite[Lemma 9]{DoLuSa:2012}) shows then that
\begin{eqnarray*}
|\{\phi_h,Q_{h,\varepsilon}^{-1} E_h\eta\}| &\le & \|\phi_h\|_0\|E_h\eta\|_0+
\Big|\{\phi_h,Q_{h,\varepsilon}^{-1}E_h\eta\}-\int_0^1 \phi_h (t)
(E_h\eta)(t)\mathrm d t\Big|\\
&\le& \|\phi_h\|_0 (\|E_h\eta\|_0+\pi h\|E_h\eta\|_{1})\le Ch ^3\|\phi_h\|_0 \|\eta\|_3,
\end{eqnarray*}
where the last inequality follows from the fact that $\|E_h\eta\|_s\le C
h^{t-s}\|\eta\|_t$ for all $t\ge s$ \cite[Theorem 8.2.1]{SaVa:2002}.
\end{proof}

For simplicity, in what follows we will write $\mathrm P\in \mathcal E(n)$ when
$\mathrm P$ is a periodic pseudodifferential operator of order $n$, i.e.,
$\mathrm P:H^s\to H^{s-n}$ is bounded for all $s$.

\begin{proposition}  \label{prop:exp2}
There exists $\mathrm P_{\mathrm k}\in \mathcal{E}(1)$  so that for all
$\eta\in H^3$ and $\psi\in H^4$, 
\begin{eqnarray*}
|\mathrm j(\phi_h,D_{h,\varepsilon}^{-1}\eta)-
\{Q_{h,\varepsilon}^{-1}\phi_h, \mathrm J\eta\}
|&\le&Ch^{3}\|\eta\|_3\|\phi_h\|_0,\quad \forall \phi_h\in
S_{h,\varepsilon},\\
|\mathrm k(\mu_h,D_h\psi)-\{\mu_h, \mathrm K
\psi\}- h^2\{\mu_h, \mathrm P_{\mathrm k}\psi\}|
&\le&Ch^{3}\|\psi\|_4\|\mu_h\|_{-1}, \quad \forall \mu_h\in
S_{h}^{-1}. 
\end{eqnarray*}
The coefficient $\mathrm P_{\mathrm k}$ and the constants in the bounds do not
depend on $\varepsilon$.
\end{proposition}

\begin{proof} We refer to \cite{DoRaSa:2008}, where similar expansions are derived. 
\end{proof}

The study of the approximation properties of $\mathrm V$ and $\mathrm W$ is
strongly influenced by the parameter $\varepsilon$. We write
\[
C_1(\varepsilon) :=\frac1{2\pi\imath}\log(\sin^2(\pi \varepsilon)) \qquad
C_2(\varepsilon):=\frac12\int_0^\varepsilon C_1(t)\mathrm d t,
\]
and note that $C_1(\pm 1/6)=0$.

\begin{proposition}    \label{prop:exp3}
There exists a smooth function $a_{\mathrm v}$ and operators $\mathrm L_{\mathrm
v}\in\mathcal{E}(1)$, $\mathrm L^1_{\mathrm w}
\in\mathcal{E}(2)$, $\mathrm L^{2a}_{\mathrm w}, \mathrm L^{2b}_{\mathrm
w}\in\mathcal{E}(3)$  such that for all $\eta\in H^3$ and $\psi\in H^4$,
\begin{eqnarray*}
&&\hspace{-0.5cm}
|\mathrm v(\mu_h,D_{h,\varepsilon}^{-1}\eta)-\{\mu_h,\mathrm V \eta\}- h C_1(\varepsilon)
\{\mu_h,a_{\mathrm v}\eta  \}\\
& & \hspace{3cm}- h^2 C_2(\varepsilon)
\{\mu_h,\mathrm L_{\mathrm v}\eta   \} |\le Ch^{3}\|\eta\|_3\|\mu_h\|_{-1},\quad \forall \mu_h\in
S^{-1}_h,\\
&&\hspace{-0.5cm}
|\mathrm w(\phi_h,D_h\psi)-\{ Q_{h,\varepsilon}^{-1} \phi_h,\mathrm W\psi\}
- h   C_1(\varepsilon)
\{ Q_{h,\varepsilon}^{-1} \phi_h ,\mathrm L_{\mathrm w}^1\psi \}\\
&&\hspace{1.6cm}- h^2
\{Q_{h,\varepsilon}^{-1}  \phi_h, (C_2(\varepsilon) \mathrm L^{2a}_{\mathrm w}+\mathrm L^{2b}_{\mathrm w})\psi 
\} \big| \le  Ch ^3\|\psi\|_4\|\phi_h\|_0,   \quad \forall \phi_h\in
S_{h,\varepsilon}.
\end{eqnarray*}
\end{proposition}
\begin{proof}
The first expansion is given in \cite[Theorem 7]{CeDoSa:2002}, while the second one is proved in 
 \cite[Proposition A.4]{DoLuSa:2012}.
\end{proof}

The key fact at this point is that by letting $\varepsilon=\pm 1/6$ 
all the expansions start at $h^2$.  This will be  crucial
since, as we will see in the next subsection, we can identify the
order of the method with the first power
of $h$ appearing in the consistency expansion. The relevance of identifying the
$h^2$ term of the asymptotic expansion of the consistency error in Propositions
\ref{prop:exp1}, \ref{prop:exp2} and \ref{prop:exp3} is related to the
possibility of moving from the norms given by the inf-sup conditions in
Proposition \ref{prop:infsup} to stronger norms when producing estimates of the
convergence error. (See Theorem \ref{the:4.7} below.)

\begin{remark}
\label{rem:4.6} 
If in Propositions \ref{prop:exp1}, \ref{prop:exp2} and \ref{prop:exp3} we only
assume that $\eta\in H^2$ and $\psi\in H^3$, and we eliminate the $h^2$ term
from the left-hand side of the bounds, then the result holds with a bound of the
form $Ch^2\|\eta\|_2$ or $Ch^2\|\psi\|_3$.
\end{remark}
\subsection{Convergence estimates}

We collect in this subsection the convergence results for the all numerical
schemes presented in this paper. 

\begin{theorem}\label{the:4.7} Assume that $k$ satisfies the hypothesis of
Proposition \ref{prop:infsup} and $\varepsilon=\pm 1/6.$ Let
$(\varphi_h,\lambda_h)\in S_{h,\varepsilon}^{-1}\times S_h$ be the pair
associated to the solution $(\boldsymbol\lambda,\boldsymbol\varphi)$ of any of
{\rm (dD01h), (dD02h), (dN01h)} or {\rm (dN02h)}. Then
\[
\|\varphi_h-D_h\varphi\|_0+\|\lambda_h-\lambda\|_{-1}\le C h^2 (\|\varphi\|_3+\|\lambda\|_2).
\]
Moreover,
\[
\max_j |\varphi_j-\beta_0(t_j)| +\max_j |h^{-1} \lambda_j-\beta_1(t_j)| \le C h^2 (\|\varphi\|_4+\|\lambda\|_4).
\]
\end{theorem}

\begin{proof}
We will only show the case (dN01h), all others being very similar. Using the
variational representation of (dN01h) in \eqref{eq:dN01h}, we can write
\[
-\smallfrac12\{\mu_h,\varphi_h\}+\mathrm k(\mu_h,\varphi_h)=\{\mu_h,-\smallfrac12 \varphi+\mathrm K\varphi\}+\{\mu_h,\mathrm V(Q_{h,\varepsilon}^{-1}\lambda-\lambda)\}\quad\forall\mu_h \in S_h^{-1}.
\]
Using now Propositions \ref{prop:exp1} (second bound), \ref{prop:exp2} (second
bound) and \ref{prop:exp3} (first bound) --see also Remark \ref{rem:4.6}-- and
\eqref{eq:Eh}, it follows that
\begin{eqnarray*}
| -\smallfrac12\{\mu_h, \varphi_h-D_h\varphi\}+\mathrm k(\mu_h,\varphi_h-D_h\varphi)|  & \le & C h^2 \|\mu_h\|_{-1} (\|\varphi\|_3+\|\lambda\|_2)\\
& &+ |\mathrm v(\mu_h,Q_{h,\varepsilon}^{-1} E_h\lambda)| \quad \forall \mu_h \in S_h^{-1}.
\end{eqnarray*}
Using \cite[Lemma 13]{CeDoSa:2002}  and the fact that
\[
\|Q_{h,\varepsilon}^{-1}\eta\|_{-1}\le C (\|\eta\|_0+ h \|\eta\|_1),
\] 
(see \cite[Lemma 6]{CeDoSa:2002})  we can prove that
\[
 |\mathrm v(\mu_h,Q_{h,\varepsilon}^{-1} E_h\lambda)|\le
 C \|\mu_h\|_{-1} ( \| E_h\lambda\|_0+ h \|E_h\lambda\|_1) \le C h^2
\|\mu_h\|_{-1}\|\lambda\|_2 \quad \forall \mu_h \in S_h^{-1}.
\]
Therefore, by Proposition \ref{prop:infsup}, the bound for
$\|\varphi_h-D_h\varphi\|_0$ follows. The bound for
\[
\|\lambda-\lambda_h\|_{-1} =\|\lambda-Q_{h,\varepsilon}^{-1}\lambda\|_{-1}\le
\|\lambda-D_{h,\varepsilon}^{-1}\lambda\|_{-1}+\|D_{h,\varepsilon}^{-1}
\lambda-Q_{h,\varepsilon}^{-1}\lambda\|_{-1}
\]
follows from \eqref{eq:Dhb} and \eqref{eq:Eh}. The uniform estimates require
including the $h^2$ term of the consistency error expansion: see \cite[Corollary
11]{CeDoSa:2002} and \cite[Theorem 6.4]{DoLuSa:2012} for very similar arguments.
\end{proof}

\begin{theorem} Assume that $k$ satisfies the hypothesis of
Proposition \ref{prop:infsup} and $\varepsilon=\pm 1/6.$ 
Let $\psi_h \in S_h$ be associated to the solution $\boldsymbol\psi$ of {\rm (iD02h)} or {\rm (iN02h)} and let $\eta_h\in S_{h,\varepsilon}^{-1}$ be associated to the solution $\boldsymbol\eta$ of {\rm (iD01h)} or {\rm (iN01h)}. Then
\[
\| D_h\psi-\psi_h\|_0\le C h^2 \|\psi\|_3 \qquad \|\eta-\eta_h\|_{-1}\le C h^2 \|\eta\|_2.
\]
\end{theorem}

\begin{proof} The proof is very similar to the one of Theorem \ref{the:4.7}. The
absence of integral operators in the right hand side makes the arguments
slightly simpler.
\end{proof}

In all cases it is possible to prove that the estimates can be transferred to
the computation of potential, with the direct representation \eqref{eq:reprh} in
the case of direct method, or the associated layer potential in the case of
indirect methods. In all cases, we can prove $|U(\mathbf z)-U_h(\mathbf z)|\le
C(\mathbf z) h^2$.

\begin{remark}
If we take $\varepsilon \neq \pm 1/6$, the methods involving $\mathrm V_h$ or $\mathrm W_h$ are of order one.
\end{remark}

\section{Experiments}

In the following experiments we consider a single elliptical obstacle with boundary
\[
\smallfrac14 (x-0.1)^2+(y-0.2)^2=1.
\]
We will check solutions in two observation points inside the ellipse $\mathbf
x_1=(0.2, 0.4)$ and $\mathbf x_2=(-0.2,-0.4)$. The examples will use more
complicated integral equations than those explained in the previous sections, in
order to put the discrete Calder\'on Calculus to a more demanding test.

\subsection{A transmission problem} Consider the coupling of the exterior
Helmholtz equation \eqref{eq:Helmholtz} with an interior equation with different
wave number
\[
\Delta V + (k/c)^2 V =0 \mbox{ in $\Omega$}
\]
(here $c>0$) and transmission conditions
\[
\gamma^+U+\beta_0=\gamma^- V, \qquad \partial_{\mathbf n}^+U+\beta_1=\alpha
\partial_{\mathbf n}^- V
\]
(with $\alpha>0$). Data are taken so that the exact solution is
\[
U(\mathbf z)= H^{(1)}_0(k|\mathbf z-\mathbf x_0|),\qquad V(\mathbf
z)=\exp(\imath (k/c) \mathbf z\cdot\mathbf d) \qquad \mathbf x_0\in \Omega,
\qquad |\mathbf d|=1.
\]
We use the symmetric formulation of Martin Costabel and Ernst Stephan
\cite{CoSt:1985} (see also \cite{LaRaSa:2009}). The main unknowns are
$\varphi^-=\gamma^- V$ and $\lambda^-=\alpha\partial_{\mathbf n}^- V$. The
system they satisfy is
\begin{equation}\label{eq:24}
\left[ \begin{array}{cc} \mathrm W_k +\alpha \mathrm W_{k/c} & \mathrm
J_k+\mathrm J_{k/c} \\[1.5ex] -\mathrm K_k-\mathrm K_{k/c}  & \mathrm
V_k+\frac1\alpha \mathrm V_{k/c}\end{array}\right]
\left[\begin{array}{c}\varphi^-\\[1.5ex] \lambda^-\end{array}\right]=
\left[\begin{array}{cc} \mathrm W_k & \smallfrac12\mathrm I+\mathrm J_k \\[1.5ex] \frac12\mathrm I-\mathrm K_k & \mathrm V_k\end{array}\right]\left[\begin{array}{c}  \beta_0 \\[1.5ex] \beta_1\end{array}\right],
\end{equation}
where we have tagged the integral operators with the corresponding wave number.
The potential representation for the interior and exterior fields is
\begin{equation}\label{eq:25}
U=-\mathrm S_k (\lambda^--\beta_1)+\mathrm D_k(\varphi^--\beta_0), \qquad V =
\alpha^{-1}\mathrm S_{k/c} \lambda^--\mathrm D_{k/c}\varphi^-.
\end{equation}
Discretization is carried out by simply substituting the elements of
\eqref{eq:24} and \eqref{eq:25} by their discrete counterparts: the data
functions are sampled with \eqref{eq:beta}, the integral operators are build
with \eqref{eq:discreteOperators} and the potentials with \eqref{eq:poth}. We
solve and tabulate the following errors:
\[
\mathrm E_h^\lambda:=\max_j | h^{-1} \lambda_j^- 
-\alpha\partial_{\mathbf n}^+ V (t_j^\varepsilon)| \qquad \mathrm E_h^\varphi:=\max_j |\varphi_j^--\gamma^- V(t_j)|
\]
\[
\mathrm E_h^V:=\max_{\ell=1,2} |V_h(\mathbf x_\ell)-V(\mathbf x_\ell)|
\]
These experiments are reported in Tables \ref{table:E1} and \ref{table:E1b}. The parameters are $k=3$, $c=2/3$ and $\alpha=3/2.$

\begin{table}[htb]
\begin{center}
\begin{tabular}{rcc}\hline
 $N$&error&e.c.r\\[0.5ex] \hline
10&4.6842$E(+000)$&\\
20&1.2470$E(+000)$&1.9093\\
40&3.7207$E(-001)$&1.7448\\
80&9.4663$E(-002)$&1.9747\\
160&2.3768$E(-002)$&1.9938\\
320&5.9518$E(-003)$&1.9976\\
640&1.4886$E(-003)$&1.9994\\
\hline
\end{tabular}\qquad
\begin{tabular}{rcc}\hline
$N$&error&e.c.r\\[0.5ex] \hline
10&5.8671$E(-001)$&\\
20&1.9979$E(-001)$&1.5542\\
40&4.9104$E(-002)$&2.0246\\
80&1.2376$E(-002)$&1.9883\\
160&3.1081$E(-003)$&1.9934\\
320&7.7699$E(-004)$&2.0001\\
640&1.9423$E(-004)$&2.0001\\
\hline
\end{tabular}
\end{center}
\caption{Errors $\mathrm E_h^\lambda$ (left column) and $\mathrm E_h^\varphi$
(right column) for the Transmission Problem in Experiment 1.}\label{table:E1}
\end{table}

\begin{table}[htb]
\begin{center}
\begin{tabular}{rcc}\hline
 $N$&error&e.c.r\\[0.5ex] \hline
10&1.8729$E(-001)$&\\
20&2.0779$E(-002)$&3.1721\\
40&4.0885$E(-003)$&2.3455\\
80&9.6559$E(-004)$&2.0821\\
160&2.4527$E(-004)$&1.9770\\
320&6.1837$E(-005)$&1.9878\\
640&1.5527$E(-005)$&1.9937\\
\hline
\end{tabular}
\end{center}
\caption{Error $E_h^V$ (potential solution $V$ at two interior observation
points) for the Transmission Problem in Experiment 1.}\label{table:E1b}
\end{table}

\subsection{Burton-Miller integral equation} Consider now the exterior Helmholtz
equation \eqref{eq:Helmholtz} with boundary condition $\gamma^+ U+\gamma
U_{\mathrm{inc}}=0$, where $\Delta U_{\mathrm{inc}}+k^2 U_{\mathrm{inc}}=0$ in a
neighborhood of the interior domain $\overline\Omega$. The well known
Burton-Miller integral equation \cite[Section 3.9]{CoKr:1983} is
\begin{equation}\label{eq:26}
\smallfrac12\xi+\mathrm J\xi+c\mathrm V \xi= \partial_{\mathbf n}U_{\mathrm{inc}}+ c \gamma U_{\mathrm{inc}}.
\end{equation}
The exterior normal derivative can be computed after solving this equation and
there are two potential representations of the solution
\begin{equation}\label{eq:27}
\lambda=\xi-\partial_{\mathbf n} U_{\mathrm{inc}} \qquad 
U=-\mathrm S \xi=-\mathrm S \lambda-\mathrm D \gamma U_{\mathrm{inc}}.
\end{equation}
The value $c=-\imath k$ is the usual choice in \eqref{eq:26}. For this value,
the equation \eqref{eq:26} is uniquely solvable independently of the frequency.
Since $\mathrm S \xi=U_{\mathrm{inc}}$ in the interior domain, we compare errors
\[
\mathrm E_h^U:=\max_{\ell=1,2} | \mathrm S_h(\mathbf x_\ell) \boldsymbol\xi- U_{\mathrm{inc}}(\mathbf x_\ell)|
\]
We also compare the density $\boldsymbol\xi$ with the solution of Problem
(dD01h) (Table \ref{tab:dd}) computing the compared error 
\[
\mathrm E_h^\xi:=\max_j | \underbrace{h^{-1} \lambda_j}_{\mbox{\footnotesize
(dD01h)}}
-\underbrace{(h^{-1}\xi_j-\partial_{\mathbf n}
U_{\mathrm{inc}}(t_j^\varepsilon))}_{\mbox{\footnotesize Burton-Miller I.E.}}|.
\]
In our numerical experiments we have taken   $U_{\mathrm{inc}}({\bf x})=
\exp(\imath k {\bf d}\cdot{\bf x})$, i.e. an acoustic plane wave,  with
direction given by the unit vector
${\bf d}=(1,1)/\sqrt2$ and  wave number $k=2$. 
The results are gathered in Table \ref{table:E2}. 

\begin{table}[htb]
\begin{center}
\begin{tabular}{rcc}\hline
 $N$&error&e.c.r\\[0.5ex] \hline
10&1.7205$E(-001)$&\\
20&3.6082$E(-002)$&2.2535\\
40&1.1990$E(-002)$&1.5894\\
80&3.7936$E(-003)$&1.6602\\
160&1.0571$E(-003)$&1.8435\\
320&2.7581$E(-004)$&1.9384\\
640&7.2185$E(-005)$&1.9339\\
\hline
\end{tabular}
\qquad
\begin{tabular}{rcc}\hline
 $N$&error&e.c.r\\[0.5ex] \hline
10&7.6790$E(+000)$&\\
20&1.8790$E(+000)$&2.0310\\
40&4.1656$E(-001)$&2.1734\\
80&8.5219$E(-002)$&2.2893\\
160&1.4703$E(-002)$&2.5351\\
320&2.2452$E(-003)$&2.7112\\
640&7.1749$E(-004)$&1.6458\\
\hline
\end{tabular}
\end{center}
\caption{Errors $\mathrm E_h^U$ (left columns) and $\mathrm E_h^\xi$ (right
columns) for the Burton-Miller integral equation in Experiment
2.}\label{table:E2}
\end{table}

\subsection{Conclusions}

We have presented a collection of compatible discretizations of the two
potentials and four boundary integral operators associated to the Helmholtz
equation on smooth parametrizable curves in the plane. We have shown discrete
stability of the discrete versions for all the operators in absence of
resonances. We have also given convergence estimates for eight integral
equations that solve the exterior Dirichlet and Neumann problems, with direct
and indirect boundary integral equations. Finally, we have tested the methods in
more complicated cases, such as systems of boundary integral equations arising
from transmission problems and combined field integral equations.

\end{document}